\documentclass[12pt,reqno,a4paper]{amsart}
\usepackage{setspace,amssymb,fourier}
\usepackage[top=20mm, bottom=20mm, left=30mm, right=30mm]{geometry}
\usepackage{ifpdf}
\ifpdf
 \usepackage[hyperindex]{hyperref}%,pagebackref
\else
 \expandafter\ifx\csname dvipdfm\endcsname\relax
 \usepackage[hypertex,hyperindex]{hyperref}%
 \else
 \usepackage[dvipdfm,hyperindex]{hyperref}%
 \fi
\fi
\theoremstyle{plain}
\newtheorem{thm}{Theorem}
\theoremstyle{remark}
\newtheorem{rem}{Remark}
\allowdisplaybreaks[4]
\DeclareMathOperator{\td}{d}

\begin{document}

\title{Six proofs for an identity of the Lah numbers}

\author[B.-N. Guo]{Bai-Ni Guo}
\address{School of Mathematics and Informatics, Henan Polytechnic University, Jiaozuo City, Henan Province, 454010, China}
\email{\href{mailto: B.-N. Guo <bai.ni.guo@gmail.com>}{bai.ni.guo@gmail.com}, \href{mailto: B.-N. Guo <bai.ni.guo@hotmail.com>}{bai.ni.guo@hotmail.com}}
\urladdr{\url{http://www.researcherid.com/rid/C-8032-2013}}

\author[F. Qi]{Feng Qi}
\address{College of Mathematics, Inner Mongolia University for Nationalities, Tongliao City, Inner Mongolia Autonomous Region, 028043, China; Department of Mathematics, College of Science, Tianjin Polytechnic University, Tianjin City, 300387, China}
\email{\href{mailto: F. Qi <qifeng618@gmail.com>}{qifeng618@gmail.com}, \href{mailto: F. Qi <qifeng618@hotmail.com>}{qifeng618@hotmail.com}, \href{mailto: F. Qi <qifeng618@qq.com>}{qifeng618@qq.com}}
\urladdr{\url{http://qifeng618.wordpress.com}}

\begin{abstract}
In the paper, utilizing respectively the induction, a generating function of the Lah numbers, the Chu-Vandermonde summation formula, an inversion formula, the Gauss hypergeometric series, and two generating functions of the Stirling numbers of the first kind, the authors collect and provide six proofs for an identity of the Lah numbers.
\end{abstract}



\thanks{Please cite this article as "Bai-Ni Guo and Feng Qi, \textit{Six proofs for an identity of the Lah numbers}, Online Journal of Analytic Combinatorics \textbf{10} (2015), 5~pages."}

\maketitle

\section{Introduction}

In combinatorics, the Lah numbers, discovered by Ivo Lah in 1955 and usually denoted by $L(n,k)$, count the number of ways a set of $n$ elements can be partitioned into $k$ nonempty linearly ordered subsets and have an explicit formula
\begin{equation}\label{a-i-k-dfn}
L(n,k)=\binom{n-1}{k-1}\frac{n!}{k!}.
\end{equation}
The Lah numbers $L(n,k)$ may also be interpreted as coefficients expressing rising factorials
\begin{equation}
(x)_n=
\begin{cases}
x(x+1)(x+2)\dotsm(x+n-1), & n\ge1\\
1, & n=0
\end{cases}
\end{equation}
in terms of falling factorials
\begin{equation}
\langle x\rangle_n=
\begin{cases}
x(x-1)(x-2)\dotsm(x-n+1), & n\ge1\\
1,& n=0.
\end{cases}
\end{equation}
For more information on the Lah numbers $L(n,k)$, please refer to, for instance, the books~\cite{Charalambides-book-2002, GKP-Concrete-Math-2nd} and the article~\cite{Lindsay-Mansour-Shattuck-JComb-2011}.
\par
In this paper, utilizing respectively the induction, a generating function of the Lah numbers, the Chu-Vandermonde summation formula, an inversion formula, the Gauss hypergeometric series, and two generating functions of the Stirling numbers of the first kind, the authors collect and provide six proofs for an identity of the Lah numbers.

\begin{thm}\label{stirling-sum-gamma-thm}
For $k\ge2$ and $n\ge0$, we have
\begin{equation}\label{Lah-No-comp-eq}
\sum_{\ell=1}^{k}(-1)^{\ell}(n+\ell)!L(k,\ell)
=
\begin{cases}
0, & 0\le n\le k-2,\\
(-1)^k\dfrac{n!(n+1)!}{(n-k+1)!}, & n\ge k-1.
\end{cases}
\end{equation}
\end{thm}

\section{Six proofs of Theorem~\ref{stirling-sum-gamma-thm}}

Now we start out to state six proofs of the identity~\eqref{Lah-No-comp-eq}.

\subsection{First proof}

For $k\ge2$ and $n\ge0$, dividing both sides of the identity~\eqref{Lah-No-comp-eq} by $k!n!$, we see that it may be rewritten equivalently as
\begin{equation}\label{equiv-iden-eq}
\sum_{\ell=1}^k(-1)^{\ell}\binom{n+\ell}{n}\binom{k-1}{\ell-1}=(-1)^{k}\binom{n+1}{k}.
\end{equation}
In~\cite[p.~169, Table~169]{GKP-Concrete-Math-2nd}, the identity~(5.24) reads that
\begin{equation}\label{Table169GKP-Concrete-Math-2nd}
\sum_{i}\binom{\ell}{m+i}\binom{s+i}{n}(-1)^i=(-1)^{\ell+m}\binom{s-m}{n-\ell}, \quad \ell\ge0.
\end{equation}
Taking $s=n$, $\ell=k-1$, and $m=-1$ in~\eqref{Table169GKP-Concrete-Math-2nd} results in~\eqref{equiv-iden-eq}.
\par
The identity~\eqref{Table169GKP-Concrete-Math-2nd} may be derived inductively, as done in~\cite[p.~170]{GKP-Concrete-Math-2nd} as follows. For $\ell=0$, all terms in~\eqref{Table169GKP-Concrete-Math-2nd} are zero except when $k=-m$, so both sides of the equation are $(-1)^m\binom{s-m}{n}$. Now suppose that the identity holds for all values less than some fixed $\ell>0$. We can use the addition formula to replace $\binom{\ell}{m+k}$ by $\binom{\ell-1}{m+k} +\binom{\ell-1}{m+k-1}$; the original sum now breaks into two sums, each of which can be evaluated by the induction hypothesis:
\begin{multline*}
\sum_k\binom{\ell-1}{m+k}\binom{s+k}{n}(-1)^k +\sum_k\binom{\ell-1}{m+k-1}\binom{s+k}{n}(-1)^k\\
=(-1)^{\ell-1+m}\binom{s-m}{n-\ell+1}+(-1)^{\ell+m}\binom{s-m+1}{n-\ell+1}.
\end{multline*}
And this simplifies to the right-hand side of~\eqref{Table169GKP-Concrete-Math-2nd}, if we apply the addition formula once again.
\par
The identity~\eqref{Table169GKP-Concrete-Math-2nd} may also be derived by reducing it to the Chu-Vandermonde summation formula
\begin{equation}\label{Chu-Vandermonde-eq}
\sum_{k}\binom{r}{m+k}\binom{s}{n-k}=\binom{r+s}{m+n}
\end{equation}
through a sequence of transformations. See~\cite[pp.~169\nobreakdash--170]{GKP-Concrete-Math-2nd}. The first proof of Theorem~\ref{stirling-sum-gamma-thm} is complete.

\subsection{Second proof}

In~\cite[p.~302, (8.40)]{Charalambides-book-2002}, it was given that
\begin{equation}\label{lah-gff-eq}
(-1)^n\langle-t\rangle_n=\sum_{k=0}^nL(n,k)\langle t\rangle_k.
\end{equation}
This is also known as the factorial generating function for the Lah numbers $L(n,k)$. Employing
\begin{align*}
\langle-m-1\rangle_km!&=(-m-1)(-m-2)\dotsm(-m-k)m!\\
&=(-1)^k(m+1)(m+2)\dotsm(m+k)m!\\
&=(-1)^k(m+k)!
\end{align*}
and replacing $t$ by $-m-1$ in the factorial generating function~\eqref{lah-gff-eq} yield
\begin{align*}
\sum_{k=0}^n(-1)^k(m+k)!L(n,k)&=(-1)^nm!(m+1)_n\\
&=(-1)^nm!(m+1)m\dotsm(m+1-n+1)\\
&=(-1)^n\frac{m!(m+1)!}{(m-n+1)!}
\end{align*}
for $m\ge n-1$. Thus, Theorem~\ref{stirling-sum-gamma-thm} is proved.

\subsection{Third proof}

Writing $\binom{k-1}{\ell-1}$ as $\binom{k-1}{k-\ell}$, we see that the identity~\eqref{equiv-iden-eq} may be obtained by extracting the coefficient of $x^k$ from both sides of the convolution product
\begin{equation*}
\frac1{(1+x)^{n+1}}(1+x)^{k-1}=\frac1{(1+x)^{n-k+2}}=\sum_{i\ge0}(-1)^i\binom{i+n-k+1}{i}x^i.
\end{equation*}
Theorem~\ref{stirling-sum-gamma-thm} is thus proved.

\subsection{Fourth proof}
The formula~\cite[p.~192, (5.48)]{GKP-Concrete-Math-2nd} reads that
\begin{equation}\label{inversion-formula}
g(k)=\sum_\ell\binom{k}{\ell}(-1)^\ell f(\ell) \Longleftrightarrow
f(k)=\sum_\ell\binom{k}{\ell}(-1)^\ell g(\ell).
\end{equation}
This dual relationship between $f$ and $g$ is called an inversion formula. Since
\begin{gather*}
\sum_{\ell=1}^{k}(-1)^{\ell}(n+\ell)!L(k,\ell)
=\sum_{\ell=1}^{k}(-1)^{\ell}(n+\ell)!\binom{k-1}{\ell-1}\frac{k!}{\ell!}\\
=\sum_{\ell=1}^{k}(-1)^{\ell}(n+\ell)!\frac{(k-1)!}{(\ell-1)!(k-\ell)!}\frac{k!}{\ell!}
=\sum_{\ell=1}^{k}\binom{k}{\ell}(-1)^{\ell}\frac{(n+\ell)!(k-1)!}{(\ell-1)!},
\end{gather*}
we may denote $f(\ell)=\frac{(n+\ell)!(k-1)!}{(\ell-1)!}$ and $f(k)=(n+k)!$. On the other hand, let
\begin{equation}
g(\ell)=(-1)^\ell\frac{(n+k-\ell)! (n+k-\ell+1)!}{(n+k-\ell+1)!},
\end{equation}
then $g(k)=(-1)^k\frac{n!(n+1)!}{(n-k+1)!}$ and
\begin{equation*}
\sum_{\ell=1}^{k}\binom{k}{\ell}(-1)^{\ell}\biggr[(-1)^\ell\frac{(n+k-\ell)! (n+k-\ell+1)!}{(n+k-\ell+1)!}\biggr]
=(n+k)!.
\end{equation*}
As a result, by the inversion formula~\eqref{inversion-formula}, we gain the identity~\eqref{Lah-No-comp-eq}.

\subsection{Fifth proof}
The generalized hypergeometric series
\begin{equation}\label{hypergeom-f}
{}_pF_q(a_1,\dotsc,a_p;b_1,\dotsc,b_q;z)=\sum_{\ell=0}^\infty\frac{(a_1)_\ell\dotsm(a_p)_\ell} {(b_1)_\ell\dotsm(b_q)_\ell}\frac{z^\ell}{\ell!}
\end{equation}
is defined for complex numbers $a_i\in\mathbb{C}$ and $b_i\in\mathbb{C}\setminus\{0,-1,-2,\dotsc\}$ and for positive integers $p,q\in\mathbb{N}$. Specially, the function ${}_2F_1(a_1,a_2;b_1,z)$ are called the Gauss hypergeometric series. See~\cite[pp.~3\nobreakdash--5]{Basic-hypergeometric-series-2nd}.
A straightforward computation reveals
\begin{gather*}
\sum_{\ell=1}^{k}(-1)^{\ell}(n+\ell)!L(k,\ell)
=\sum_{\ell=0}^\infty(-1)^{\ell}(n+\ell)!L(k,\ell)\\
=-k! (n+1)!{}_2F_1(1-k,n+2;2;1).
\end{gather*}
Applying the Chu-Vandermonde summation formula~\eqref{Chu-Vandermonde-eq} to the Gauss hypergeometric series ${}_2F_1(1-k,n+2;2;1)$ results in the identity~\eqref{Lah-No-comp-eq}.

\subsection{Sixth proof}

After collecting the above proofs, we now provide a new proof of the identity~\eqref{Lah-No-comp-eq}. Although this proof is not simple, we believe that it is still worthwhile to announce it, because it contains some novel idea and comprehensive understanding on the Lah numbers.
\par
It is easy to verify that
\begin{equation}\label{frac1-ln-x+1-int}
\frac1{\ln(1+t)}=\int_0^\infty\frac1{(1+t)^u}\td u
\end{equation}
for $t>0$.
In~\cite[Lemma~2]{Liu-Qi-Ding-2010-JIS} and~\cite[Section~2]{Filomat-36-73-1.tex}, it was obtained inductively that
\begin{equation}\label{Liu-Qi-Ding-2010-JIS-log-der}
\biggl[\frac1{\ln(1+t)}\biggr]^{(m)} =\frac1{(1+t)^m}\sum_{k=0}^m(-1)^k\frac{k!s(m,k)}{[\ln(1+t)]^{k+1}}, \quad m\ge0,
\end{equation}
where $s(m,k)$ denote the Stirling numbers of the first kind and may be generated by
\begin{equation}\label{gen-funct-3}
\frac{[\ln(1+t)]^k}{k!}=\sum_{n=k}^\infty s(n,k)\frac{t^n}{n!},\quad |t|<1.
\end{equation}
Utilizing the integral representation~\eqref{frac1-ln-x+1-int} in the left hand side of equation~\eqref{Liu-Qi-Ding-2010-JIS-log-der} and simplifying give
\begin{equation}\label{gamma-stirling}
\int_0^\infty\frac{\Gamma(u+m)}{\Gamma(u)}\frac1{(1+t)^{u}}\td u =(-1)^{m}\sum_{i=0}^m(-1)^{i}\frac{i!s(m,i)}{[\ln(1+t)]^{i+1}},
\end{equation}
where $\Gamma(z)$ is the classical Euler gamma function which may be defined by
\begin{equation}\label{gamma-dfn}
\Gamma(z)=\int^\infty_0u^{z-1} e^{-u}\td u, \quad \Re(z)>0.
\end{equation}
Substituting $t$ for $\frac1{\ln(1+t)}$ in~\eqref{gamma-stirling} brings out
\begin{equation}\label{gamma-stir-t}
\int_0^\infty\frac{\Gamma(u+m)}{\Gamma(u)}e^{-u/t}\td u =(-1)^{m}\sum_{i=0}^m(-1)^{i}i!s(m,i)t^{i+1}.
\end{equation}
Differentiating $1\le k\le m+1$ times with respect to $t$ on both sides of~\eqref{gamma-stir-t} generates
\begin{equation}\label{1stirl-gamma-int}
\int_0^\infty\frac{\Gamma(u+m)}{\Gamma(u)}\bigl(e^{-u/t}\bigr)^{(k)}\td u =(-1)^{m}\sum_{i=k-1}^m(-1)^{i}\frac{i!(i+1)!}{(i-k+1)!}s(m,i)t^{i-k+1}.
\end{equation}
\par
In~\cite[Theorem~2.2]{exp-reciprocal-cm-IJOPCM.tex}, it was obtained and applied that
\begin{equation}\label{g(t)-derivative}
\bigl(e^{-1/t}\bigr)^{(i)}=\frac1{e^{1/t}t^{2i}} \sum_{k=0}^{i-1}(-1)^kL(i,i-k){t^{k}}, \quad i\in\mathbb{N}.
\end{equation}
See also~\cite{DMST-MM-2013-Exp}.
Consequently, it follows that
\begin{equation*}
\bigl(e^{-u/t}\bigr)^{(k)}=\frac{u^{k}}{e^{u/t}t^{2k}} \sum_{\ell=0}^{k-1}(-1)^\ell\frac{L(k,k-\ell)}{u^\ell}{t^{\ell}}.
\end{equation*}
Substituting this formula into~\eqref{1stirl-gamma-int} and making use of another generating function of the Stirling numbers of the first kind $s(n,k)$
\begin{equation}\label{gen-funct-1}
\frac{\Gamma(x+n)}{\Gamma(x)}
=\prod_{k=0}^{n-1}(x+k)
=\sum_{k=0}^n(-1)^{n-k}s(n,k)x^k
\end{equation}
reveals
\begin{align*}
&\quad(-1)^{m}\sum_{i=k-1}^m(-1)^{i}\frac{i!(i+1)!}{(i-k+1)!}s(m,i)t^{i-k+1}\\
&=\frac1{t^{2k}}\sum_{\ell=0}^{k-1}(-1)^{\ell}L(k,k-\ell)t^{\ell} \int_0^\infty\frac{\Gamma(u+m)}{\Gamma(u)}u^{k-\ell}e^{-u/t}\td u\\
&=\frac1{t^{2k}}\sum_{\ell=0}^{k-1}(-1)^{\ell}L(k,k-\ell)t^{\ell}
\sum_{i=0}^m(-1)^{m-i}s(m,i)\int_0^\infty u^{i+k-\ell}e^{-u/t}\td u\\
&=\sum_{\ell=0}^{k-1}\sum_{i=0}^m(-1)^{\ell}L(k,k-\ell)
(i+k-\ell)!(-1)^{m-i}s(m,i)t^{i-k+1}\\
&=(-1)^{m}\sum_{i=0}^m(-1)^{i}\Biggl[\sum_{\ell=0}^{k-1}(-1)^{\ell}(i+k-\ell)! L(k,k-\ell)\Biggr]s(m,i)t^{i-k+1}.
\end{align*}
Equating coefficients of the factors $t^{i-k+1}$ in the above equation produces
\begin{equation*}
\sum_{\ell=0}^{k-1}(-1)^{\ell}(n+k-\ell)!L(k,k-\ell)
=
\begin{cases}
0, & 0\le n\le k-2\\
\dfrac{n!(n+1)!}{(n-k+1)!}, & n\ge k-1
\end{cases}
\end{equation*}
which may be rearranged as~\eqref{Lah-No-comp-eq}. Theorem~\ref{stirling-sum-gamma-thm} is thus proved.

\begin{rem}
It is interesting that, the Stirling numbers of the first kind $s(n,k)$ and two of their generating functions~\eqref{gen-funct-3} and~\eqref{gen-funct-1} are utilized in the proof of Theorem~\ref{stirling-sum-gamma-thm}, but they do not appear in the final result~\eqref{Lah-No-comp-eq}.
\end{rem}

\begin{rem}
This paper is an extended version of the preprint~\cite{Lah-No-Identity.tex}.
\end{rem}

\subsection*{Acknowledgements}
The authors would like to express many thanks to Professor Dr Yi Wang at Dalian University of Technology in China, to Professor Dr Christian Krattenthaler at the Fakult\"at f\"ur Mathematik of the Universit\"at Wien in Austria, and to two anonymous mathematicians for their discussions about, comments on, and contributions to the original version of this paper.

\end{document}